\documentclass[10pt,draft]{article}
\usepackage{amsmath,amsfonts,amssymb,amsthm,amsxtra,bookman,latexsym,verbatim}
\usepackage[mathscr]{eucal}
\usepackage[all]{xy}

\theoremstyle{plain}

\theoremstyle{definition}

\theoremstyle{remark}

\newcommand{\N}{\mathbb{N}}
\newcommand{\Q}{\mathbb{Q}}
\newcommand{\R}{\mathbb{R}}
\newcommand{\Z}{\mathbb{Z}}

\renewcommand{\P}{\mathcal{P}}
\renewcommand{\S}{\mathcal{S}}

\renewcommand{\a}{\mathfrak{a}}

\newcommand{\p}{\mathfrak{p}}
\newcommand{\q}{\mathfrak{q}}

\DeclareMathOperator{\Bl}{Bl} \DeclareMathOperator{\Cl}{Cl} \DeclareMathOperator{\Div}{Div}   \DeclareMathOperator{\Hom}{Hom}   \DeclareMathOperator{\ord}{ord}     \DeclareMathOperator{\Spec}{Spec}   

\begin{document}

\title{Comments on toric varieties}
\author{Howard~M Thompson \thanks{hmthomps@umich.edu}}
\maketitle

\begin{abstract}
   Here are few notes on not necessarily normal toric varieties and resolution by toric blow-up. These notes are
   independent of, but in the same spirit as the earlier preprint \cite{hT03a}. That is, they focus on the fact that toric
   varieties are locally given by monoid algebras.
\end{abstract}


\section{Not necessarily normal toric varieties}

Fix a base field $k$. We start by recalling the definition of a fan. Let $N\cong\Z^d$ be a lattice, that is, a finitely generated free Abelian group. A (strongly convex rational polyhedral) cone, $\sigma$, in $N_\R=N\bigotimes_\Z\R$ is a set consisting of all nonnegative linear combinations of some fixed finite set of vectors in the lattice, 
\[
   \sigma=\R_{\geq0}v_1+\cdots\R_{\geq0}v_r,\qquad v_1,\ldots,v_r\in N
\]
that contains no line. Here we identify $N$ with its image, $\{n\otimes1\mid n\in N\}$, in $N_\R$. 

Let $M=\Hom(N,\Z)$ be the dual lattice to $N$, identify $M$ with its image in $M_\R$, identify $M_\R$ with the dual space to $N_\R$, and let $\langle,\rangle$ be the dual pairing.

We say a $d-1$-dimensional subspace, $H$, of $N_\R$ is a supporting hyperplane of $\sigma$ if there exists a vector $u\in M_\R$ such that $H=\{v\mid\langle u,v\rangle=0\}$ and $\sigma\subset\{v\mid\langle u,v\rangle\geq0\}$. A face of a cone $\sigma$ is a subset of the form $H\cap\sigma$ where $H$ is a supporting hyperplane of $\sigma$.

A fan, $\Delta$ is a finite collection of cones that is closed under taking faces such that the intersection of any two cones in $\sigma$ is a face of each.

To each cone $\sigma$, we associate: (1) a finitely generated submonoid $S_{\sigma}=\sigma\spcheck\cap M$ of $M$, where $\sigma\spcheck=\{u\in M_\R\mid\langle u,v\rangle\geq0,\,\forall v\in\sigma\}$; (2) the finitely generated $k$-algebra $k[S_{\sigma}]$; and, (3) the affine $k$-variety $U_{\sigma}=\Spec k[S_{\sigma}]$. The (affine) toric variety associated to $\sigma$ is $U_{\sigma}$. If $\tau$ is a face of $\sigma$, $k[S_{\tau}]$ is a localization of $k[S_{\sigma}]$ and $U_{\tau}$ is an open affine subset of $U_{\sigma}$. Using these identifications, we associate an algebraic variety to a fan $\Delta$. We call this variety, $X_{\Delta}$, the toric variety associated to $\Delta$.

We say a submonoid $S\subseteq M$ is saturated if $S=\R_{\geq0}S\cap\Z S$. That is, a saturated monoid is the intersection of the lattice it generates with the cone it generates in the real vector space it generates. The monoids $S_{\sigma}$ are saturated. For a finitely generated submonoid of $M$, we call the monoid $S^{sat}=\R_{\geq0}S\cap\Z S$ the saturation of $S$. In fact, $S^{sat}=\Q_{\geq0}S\cap\Z S=\{u\in M\mid nu\in S\text{ for some positive integer }n$. If $S$ is a finitely generated submonoid of $M$, then $S^{sat}$ is a finitely generated saturated submonoid of $M$. Hochster~\cite{mH72} proved the monoid algebra of a finitely generated saturated submonoid of $M$ is integrally closed. Evidently, $k[S^{sat}]$ is integral over $k[S]$. So, $k[S^{sat}]$ is integral the integral closure of $k[S]$.

In order to give a not necessarily normal version of toric varieties, we will abandon this description in terms of fans. More specifically, the duality in the step $\sigma\rightsquigarrow S_{\sigma}$ forces the normality of the scheme $X_{\Delta}$. Our approach will be to characterize the set of monoids $\{S_{\sigma}\mid\sigma\in\Delta\}$ and then consider collections of monoids that satisfy all the conditions of our characterization except that of saturation.

First, note that such a set $\{S_{\sigma}\mid\sigma\in\Delta\}$ consists of finitely generated saturated submonoids $S\subseteq M$ such that $\Z S=M$. We will also rely on the following two facts. If $\sigma$ and $\tau$ are any two cones in $N$, then $(\sigma\cap\tau)\spcheck=\sigma\spcheck+\tau\spcheck$ (see Ewald~\cite[V.2.2~Lemma]{gE96}). And, if $S_{\tau}$ is a localization of $S_{\sigma}$, then $\tau$ is a face of $\sigma$. To see this, take $u\in M$ such that $S_{\tau}=S_{\sigma}+\N u$. I claim $\tau=\{v\in\sigma\mid\langle u,v\rangle=0\}$.

In light of these facts, here is the promised characterization: Let $N$ be a lattice with dual lattice $M$ and let $\S$ be a finite collection of finitely generated saturated submonoids of $M$ such that each $S\in\S$ generates $M=\Z S$. Then, there exists a fan $\Delta$ in $N$ such that $\S=\{S_{\sigma}\mid\sigma\in\Delta\}$ if and only if $\S$ is closed under localization and the sum of any two elements of $\S$ is a localization of each. 

Let $\S$ be a finite collection of finitely generated submonoids of a lattice $M$ that is closed under localization such that the sum of any two elements of $\S$ is a localization of each and such that each element of $\S$ generates $M$. The (not necessarily normal) toric variety associated to $\S$ is obtained in the same manner as in the normal case. Such a collection yields a (generalized) fan $\S^{top}$ in the sense of Thompson~\cite{hT03a} by gluing the spectra of the monoids using the same prescription. If all the monoids in $\S$ are saturated, then this topological space is just the orbit space of the toric variety equipped with a sheaf of monoids. Henceforth, we will treat such collections as if they were the fans and we will write $X_{\Delta}$ for the (not necessarily normal) toric variety associated to the collection of monoids $\Delta$. We should note that the schemes formed this way really are varieties. The normalization of $X_{\Delta}$ is the normal toric variety associated to the fan in $\Hom(M,\Z)$ obtained by taking the collection $\{(\R_{\geq0}S)\spcheck\mid S\in\Delta\}$. In particular, since the normalization of $X_{\Delta}$ is separated, so is $X_{\Delta}$. And, $k[S]$ is a domain for each $S\in\Delta$ since it is a subring of the domain $k[M]\cong k[t_1,t_1^{-1},t_2,t_2^{-1},\ldots,t_d,t_d^{-1}]$.

In addition, the normalization map is a blow-up. Let $S$ is a finitely generated submonoid of $\Z S=M$. Suppose $s,s'\in S$ are such that $s=s'-s''\in S^{sat}$. Write $t^{s'}$ (resp. $t^{s''}$) for the image of $s'$ (resp. $s''$) in $k[S]$ and consider the affine patches of the blow-up $\Bl_{(t^{s'},t^{s''})}(\Spec k[S])$ associated to $t^{s'}$ and $t^{s''}$. Since $ns\in S$ for some positive integer $n$, the patch obtained by making $t^{s''}$ the principle generator is $\Spec k[S+\N s]$. And, $S\subseteq S+\N s\subseteq S^{sat}$. On the other hand, the patch given by making the principal generator of the ideal sheaf $t^{s'}$ is $\Spec k[S+\N s]_{t^{ns}}$, an open subset of the other patch. In particular, $s=ns+(n-1)(-s)\in S+\N(-s)$ so $S+\N s\subseteq S+\N(-s)$, $ns\in S$ becomes invertible when one adjoins $-s$ to $S$, and $-s=(n-1)s+(-ns)$ so $S+\N(-s)\subseteq S+\N s+\N(-ns)$. Since $S^{sat}$ is finitely generated and a composition of blow-ups is a blow-up, the normalization is a blow-up. In fact, fix a finite generating set $s_1,s_2,\ldots,s_m$ for $S^{sat}$ and a choice of pairs of elements $s'_i.s''_i\in S$ such that $s_i=s'_i-s''_i$ for each $i$. Now, let $I=\prod_{i=1}^m(t^{s'_i},t^{s''_i})$. Then, $\Spec k[S^{sat}]\cong\Bl_I(\Spec k[S])$. So, the normalization is also toric map in the sense of Thompson~\cite{hT03a}.

\section{Toric blow-ups and the toric variety associated to a lattice polyhedron}

For the rest of this paper, all toric varieties are normal unless stated otherwise.

Let $M$ be a lattice and let $\P$ be a full dimensional polyhedron in $M_\R$. That is, $\P$ is an intersection of finitely many half-spaces with nonempty interior. We will say $\P$ is a lattice polyhedron if for all $d$ every $d$-face of $\P$ contains $d+1$ affinely independent lattice points. Let $\P$ be a lattice polyhedron. By replacing $\P$ with some positive integer multiple of $\P$ if necessary, we may assume every face of $\P$ has a lattice point in its relative interior. We will now describe a collection of submonoids of $M$ associated to $\P$ in such a way that when $\P$ is a polytope (that is, when $\P$ is bounded) the toric variety obtained this way is (abstractly) isomorphic to the projective toric variety that is traditionally associated to this polytope. To this end, for each face $F$ of $\P$, fix a lattice point in its relative interior $u_F$. We associate the monoid $S_F=\R_{\geq0}(\P-u_F)\cap M$ where $\P-u_F=\{u-u_F\mid u\in\P\}$ to the face $F$ and the set $\{S_F\mid F\text{ is a face of }\P\}$ to $\P$. The toric variety obtained this way is quasi-projective. To see this when $\P$ is unbounded, further intersect $\P$ with a half-space in such a way as to obtain a polytope that has facets parallel to each facet of $\P$. I claim the toric variety associated to $\P$ is the open subvariety obtained from the toric variety of the polytope by removing the divisor corresponding to the facet contained in the hyplane bounding the new half-space.

We will now give a local description of toric blow-up for toric varieties. Let $S$ be a finitely generated saturated submonoid of the lattice $M=\Z S$ and let $\a$ be an integrally closed ideal of $S$. That is, $\a\subseteq S$, $\a+S=\a$, and $\a$ is the intersection of the convex hull of $\a$ and $M$ in $M_\R$. Let $I\subseteq k[S]$ be the ideal generated by $\{t^s\mid s\in\a\}$. Since the convex hull of $\a$ is a lattice polyhedron, we have two ways to associate a toric variety to $\a$. We could take the toric variety associated to the convex hull or we could take the blow-up $\Bl_I(\Spec k[S])$. These two toric varieties are isomorphic. 

Here is a sketch of the proof: We may replace $I$ with a power of $I^n$ without changing the blow-up and we may replace $\a$ with $n\a$ without changing the toric variety assicated to the convex hull. So, we may assume the relative interior of each face of the convex hull contains a lattice point by making such simultaneous replacements with $n$ large enough. Fix a generating set for $\a$ and notice that if $s$ is one of the generators, then the affine patch of the blow-up where $t^s$ is principle is isomorphic to $\Spec k[S_F]$ where $F$ is the unique face of the convex hull such that $s$ is in its relative interior. In other words, this patch is isomophic to $\Spec k[S_F]$ where $F$ is the smallest face of the convex hull containing $s$.

In particular, the faces of the convex hull of any integrally closed ideal $\a\subseteq S$ are in inclusion preserving bijection with the torus invariant pieces of $\Bl_I(\Spec k[S])$ where $I=(t^s)_{s\in\a}$.

\section{Simplicialization of non-simplicial normal toric varieties}

Let $X=X_{\Delta}$ be a toric variety. We will consider the cokernel of the standard map from the Picard group of $X$ to the (Weil) divisor class group of $X$. Or equivalently, we will consider the cokernel of the standard map from the torus invariant Cartier divisors to the torus invariant Weil divisors:
\[
   \xymatrix@C=0.5cm{
      0 \ar[r] & \Div_TX \ar[rr] && \bigoplus_{i=1}^r\Z\cdot D_i \ar[rr] && G \ar[r] & 0 }
\]
$G$ is a finitely generated Abelian group. This group is finite if and only if $X$ is simplicial. And, it is trivial if and only if $X$ is smooth. In other words, a toric variety is simplicial if and only if every Weil divisor is $\Q$-Cartier. Furthermore, a toric variety is smooth if and only if every Weil divisor is Cartier.

To see $X$ is simplicial when $G$ is finite, we work locally: Let $S\neq M$ be a finitely generated saturated submonoid of a lattice $M=\Z S$ and let $\p$ be an $S$-graded height one prime of $k[S]$. Here $\a=\{s\in S\mid t^s\in\p$ is a prime ideal of $S$. That is, $s\in\a$ or $s'\in\a$ whenever $s+s'\in\a$ and $\a+S=\a$. The complement of $\a$ in $S$ generates a supporting hyperplane $H_{\p}$ of $\R_{\geq0}S$ such that $F_{\p}=H_{\p}\cap\R_{\geq0}S$ is a facet (maximal proper face) of $\R_{\geq0}S$. For any $u\in M$, $\ord_{\p}(t^u)$ is, up to sign, the lattice distance from $s$ to $H$. A positive integer multiple $mD$ of the divisor $D$ corresponding to $\p$ consists of the $m$th symbolic power of $\a$, $\a^{(m)}=\{s\in S\mid\ord_{\p}(t^s)\geq m\}$. Think of this as the lattice points in the cone over $S$ that lie above the hyperplane parallel to $H_{\p}$ at lattice height $m$. If the image of $D$ has finite order $m$ in $G$, then $\a^{(m)}$ is a principal ideal. In this case, the principal generator of $\a^{(m)}$ must lie in the facet $F_{\q}$ for every $S$-graded height one prime $\q\neq\p$ of $k[S]$. This is due to the fact that every such facet contains lattice ponts arbitrarily far away from $H_{\p}$ and both $s$ \& $s'$ must lie in a face $F$ of the cone over $S$ whenever $s+s'\in F$ since the faces of this cone are exactly the complements of the prime ideals of $S$. In particular, if $G$ is finite, the intersection of all but one of the facets of the cone over $S$ contains a nonzero lattice point for every $S\in\Delta$. This forces each of these cones to be simplicial. The converse is a standard fact.

In the special case when $X=\Spec k[S]$, we write $\Cl(S)$ for the cokernel because it is the divisor class group of $S$. It is straightforward to see that if $D$ is a torus invariant Weil divisor whose image in $G$ has finite order, then this order is the least common multiple of the orders of the primes $\p$ corresponding to $D$ in each $\Cl(S)$ such that the point corresponding to $S$ in the generalized fan the orbit space) lies on the image of $D$. In particular, if $X_{\sigma}$ is the toric variety associated to a simplicial cone $\sigma$ in $N$, then the order of our group is the multiplicity of $\sigma$. The claim, ``A toric variety is smooth if and only if every Weil divisor is Cartier'' is an easy consequence of this fact.

Let $X=\Spec k[S]$ be a non-simplicial toric variety, let $\p$ is a height one prime whose image in $\Cl(S)$ has infinite order, let $m$ be the least common multiple of the heights of the first lattice points on each one-dimensional face of the cone over $S$ not contained in the facet $F_{\p}$, and let $\a$ be as before. In this case, the convex hull of $\a^{(m)}$ is the half-space above the hyperplane parallel to $H_{\p}$ at a lattice height $m$ above $H_{\p}$ and the cone over $S$. So, there is a one-to-one correspondence between the torus invariant divisors of $X$ and those of the blow-up of the $m$th symbolic power of $\p$. Therefore, this blow-up $\widetilde{X}\to X$ makes $X$ more simplicial without introducing new torus invariant Weil divisors because the rank of $G$ goes down. More generally, if $D$ is a torus invariant Weil divisor on a toric variety whose image in $G$ has infinite order, for each $S$ on $D$ we have this number $m$ as in the affine case. Let $m'$ be the least common multiple of these numbers. When we blow-up $m'D$ we get a toric variety that is more simplicial without introducing any new torus invariant Weil divisors. Repeatedly doing this simplicializes $X$ without introducing new invariant Weil divisors. A study of how the geometry a non-simplicial toric variety is reflected in the finitely many simplicial toric varieties obtained this way might prove interesting.

It is difficult to find the following fact in the literature: If $X$ is a toric variety, then there exists a toric resolution of singularities $\pi:\widetilde{X}\to X$ such that $\pi$ is a projective morphism. I discovered this fact through the considerations above. In hindsight though, this what happens if one uses only steller subdivision in the resolution in standard description given by fans. Our simplicializations exactly correspond to taking rays in the fan in $N$ that are contained in nonsimplicial cones and one at a time taking steller subdivisions along them. The standard way to resolve simplicial toric varieties is by steller subdivision. We have nothing new to add other than specifying which ideals are being blown-up.

\providecommand{\bysame}{\leavevmode\hbox to3em{\hrulefill}\thinspace}
\providecommand{\MR}{\relax\ifhmode\unskip\space\fi MR }
\providecommand{\MRhref}[2]{%
  \href{http://www.ams.org/mathscinet-getitem?mr=#1}{#2}
}
\providecommand{\href}[2]{#2}

\end{document}